\documentclass[a4paper,12pt,draft]{amsart}
\usepackage{amsmath, amssymb,braket}
\newcommand{\url}[1]{{\tt #1}}

\newcommand{\PPP}{{\bf{P}}}
\newcommand{\QQQ}{{\bf{Q}}}
\newcommand{\PPBBTT}{{\bf{PBT}}}

\newcommand{\ZZ}{\mathbb{Z}}
\newcommand{\NN}{\mathbb{N}}

\newcommand{\nodeminus}{\circleddash}
\newcommand{\stripminus}{\ominus}
\newcommand{\Trees}{\mathbb{T}}
\newcommand{\numof}[1]{|#1|}
\newcommand{\denotes}{:=}

\newcommand{\SU}[3][{}]{s^{#2,#3}_{#1}}
\newcommand{\SD}[3][{}]{s_{#2,#3}^{#1}}

\newcommand{\defit}[1]{{\em #1}}

\theoremstyle{plain}
\newtheorem{thm}{Theorem}[section]
\newtheorem{lemma}[thm]{Lemma}
\newtheorem{cor}[thm]{Corollary}
\newtheorem{theorem}[thm]{Theorem}

\theoremstyle{remark}
\newtheorem{remark}[thm]{Remark}

\newtheorem{example}[thm]{Example}

\theoremstyle{definition}
\newtheorem{definition}[thm]{Definition}

\allowdisplaybreaks[4]

\newcommand{\examplei}{
For  
\begin{align*}
T=&\raisebox{-84\unitlength}{\begin{picture}(58, 84)(5, -36)
\put(0,0){\line(1,2){6}}
\put(6,12){\line(1,2){6}}
\put(18,12){\line(-1,2){6}}
\put(12,24){\line(1,1){12}}
\put(30,-12){\line(-1,2){6}}
\put(24,0){\line(1,2){6}}
\put(36,0){\line(-1,2){6}}
\put(30,12){\line(1,2){6}}
\put(36,24){\line(-1,1){12}}
\put(24,36){\line(2,1){24}}
\put(48,12){\line(-1,1){12}}
\put(42,0){\line(1,2){6}}
\put(36,-12){\line(1,2){6}}
\put(30,-24){\line(1,2){6}}
\put(36,-36){\line(-1,2){6}}
\put(0,0){\circle{3}}
\put(6,12){\circle{3}}
\put(18,12){\circle{3}}
\put(12,24){\circle{3}}
\put(30,-12){\circle{3}}
\put(24,0){\circle{3}}
\put(36,0){\circle{3}}
\put(30,12){\circle{3}}
\put(36,24){\circle{3}}
\put(24,36){\circle{3}}
\put(30,-24){\circle{3}}
\put(48,12){\circle*{3}}
\put(42,0){\circle*{3}}
\put(36,-12){\circle{3}}
\put(36,-36){\circle{3}}
\put(48,48){\circle*{3}}
\put(48,48){\circle*{2}}
\put(48,12){\circle*{2}}
\put(42,0){\circle*{2}}
\put(30,36){\vector(2,1){16}}
\put(40,-12){\vector(1,2){10}}
 \end{picture}}
,& T \stripminus r_{T,3}=&%
\raisebox{-48\unitlength}{\begin{picture}(48, 60)(0, -12)
\put(0,0){\line(1,2){6}}
\put(6,12){\line(1,2){6}}
\put(18,12){\line(-1,2){6}}
\put(12,24){\line(1,1){12}}
\put(30,-12){\line(-1,2){6}}
\put(24,0){\line(1,2){6}}
\put(36,0){\line(-1,2){6}}
\put(30,12){\line(1,2){6}}
\put(36,24){\line(-1,1){12}}
\put(48,12){\line(-1,1){12}}
\put(42,0){\line(1,2){6}}
\put(48,-12){\line(-1,2){6}}
\put(0,0){\circle{3}}
\put(6,12){\circle{3}}
\put(18,12){\circle{3}}
\put(12,24){\circle{3}}
\put(30,-12){\circle{3}}
\put(24,0){\circle{3}}
\put(36,0){\circle{3}}
\put(30,12){\circle{3}}
\put(36,24){\circle{3}}
\put(24,36){\circle{3}}
\put(42,0){\circle{3}}
\put(48,12){\circle{3}}
\put(48,-12){\circle{3}}
 \end{picture}}
\end{align*}
since 
\begin{gather*}
\raisebox{-84\unitlength}{\begin{picture}(58, 84)(5, -36)
\put(0,0){\line(1,2){6}}
\put(6,12){\line(1,2){6}}
\put(18,12){\line(-1,2){6}}
\put(12,24){\line(1,1){12}}
\put(30,-12){\line(-1,2){6}}
\put(24,0){\line(1,2){6}}
\put(36,0){\line(-1,2){6}}
\put(30,12){\line(1,2){6}}
\put(36,24){\line(-1,1){12}}
\put(24,36){\line(2,1){24}}
\put(48,12){\line(-1,1){12}}
\put(42,0){\line(1,2){6}}
\put(36,-12){\line(1,2){6}}
\put(30,-24){\line(1,2){6}}
\put(36,-36){\line(-1,2){6}}
\put(0,0){\circle{3}}
\put(6,12){\circle{3}}
\put(18,12){\circle{3}}
\put(12,24){\circle{3}}
\put(30,-12){\circle{3}}
\put(24,0){\circle{3}}
\put(36,0){\circle{3}}
\put(30,12){\circle{3}}
\put(36,24){\circle{3}}
\put(24,36){\circle{3}}
\put(30,-24){\circle{3}}
\put(48,12){\circle*{3}}
\put(42,0){\circle*{3}}
\put(36,-12){\circle{3}}
\put(36,-36){\circle{3}}
\put(48,48){\circle*{3}}
\put(48,48){\circle*{2}}
\put(48,12){\circle*{2} }
\put(42,0){\circle*{2}}
\put(40,-2){\dashbox(4,4){}}
\put(40,-12){\vector(1,2){5}}
 \end{picture}}\to
\raisebox{-84\unitlength}{\begin{picture}(58, 84)(5, -36)
\put(0,0){\line(1,2){6}}
\put(6,12){\line(1,2){6}}
\put(18,12){\line(-1,2){6}}
\put(12,24){\line(1,1){12}}
\put(30,-12){\line(-1,2){6}}
\put(24,0){\line(1,2){6}}
\put(36,0){\line(-1,2){6}}
\put(30,12){\line(1,2){6}}
\put(36,24){\line(-1,1){12}}
\put(24,36){\line(2,1){24}}
\put(48,12){\line(-1,1){12}}
\put(42,0){\line(1,2){6}}
\put(36,-12){\line(1,2){6}}
\put(42,-24){\line(-1,2){6}}
\put(0,0){\circle{3}}
\put(6,12){\circle{3}}
\put(18,12){\circle{3}}
\put(12,24){\circle{3}}
\put(30,-12){\circle{3}}
\put(24,0){\circle{3}}
\put(36,0){\circle{3}}
\put(30,12){\circle{3}}
\put(36,24){\circle{3}}
\put(24,36){\circle{3}}
\put(48,12){\circle*{3}}
\put(46,10){\dashbox(4,4){}}
\put(42,0){\circle{3}}
\put(36,-12){\circle{3}}
\put(42,-24){\circle{3}}
\put(48,48){\circle*{3}}
\put(46,0){\vector(1,2){5}}
 \end{picture}}\to
\raisebox{-84\unitlength}{\begin{picture}(58, 84)(5, -36)
\put(0,0){\line(1,2){6}}
\put(6,12){\line(1,2){6}}
\put(18,12){\line(-1,2){6}}
\put(12,24){\line(1,1){12}}
\put(30,-12){\line(-1,2){6}}
\put(24,0){\line(1,2){6}}
\put(36,0){\line(-1,2){6}}
\put(30,12){\line(1,2){6}}
\put(36,24){\line(-1,1){12}}
\put(24,36){\line(2,1){24}}
\put(48,12){\line(-1,1){12}}
\put(42,0){\line(1,2){6}}
\put(48,-12){\line(-1,2){6}}
\put(0,0){\circle{3}}
\put(6,12){\circle{3}}
\put(18,12){\circle{3}}
\put(12,24){\circle{3}}
\put(24,0){\circle{3}}
\put(36,0){\circle{3}}
\put(30,12){\circle{3}}
\put(36,24){\circle{3}}
\put(24,36){\circle{3}}
\put(48,12){\circle{3}}
\put(42,0){\circle{3}}
\put(48,-12){\circle{3}}
\put(48,48){\circle*{3}}
\put(46,46){\dashbox(4,4){}}
\put(30,36){\vector(2,1){16}}
 \end{picture}}\to
\raisebox{-48\unitlength}{\begin{picture}(48, 60)(0, -12)
\put(0,0){\line(1,2){6}}
\put(6,12){\line(1,2){6}}
\put(18,12){\line(-1,2){6}}
\put(12,24){\line(1,1){12}}
\put(30,-12){\line(-1,2){6}}
\put(24,0){\line(1,2){6}}
\put(36,0){\line(-1,2){6}}
\put(30,12){\line(1,2){6}}
\put(36,24){\line(-1,1){12}}
\put(48,12){\line(-1,1){12}}
\put(42,0){\line(1,2){6}}
\put(48,-12){\line(-1,2){6}}
\put(0,0){\circle{3}}
\put(6,12){\circle{3}}
\put(18,12){\circle{3}}
\put(12,24){\circle{3}}
\put(30,-12){\circle{3}}
\put(24,0){\circle{3}}
\put(36,0){\circle{3}}
\put(30,12){\circle{3}}
\put(36,24){\circle{3}}
\put(24,36){\circle{3}}
\put(42,0){\circle{3}}
\put(48,12){\circle{3}}
\put(48,-12){\circle{3}}
 \end{picture}}.
\end{gather*}
The natural inclusion $\nu_{T,3}$ maps the nodes $\circ$ in
$T\stripminus r_{T,3}$
to the former nodes $\circ$ in $T$.
}
\newcommand{\exampleii}{
For example, for $w=1221$ and
\begin{align*}
T=&
\raisebox{-84\unitlength}{\begin{picture}(58, 84)(-5, -36)
\put(0,0){\line(1,2){6}}
\put(6,12){\line(1,2){6}}
\put(18,12){\line(-1,2){6}}
\put(12,24){\line(1,1){12}}
\put(30,-12){\line(-1,2){6}}
\put(24,0){\line(1,2){6}}
\put(36,0){\line(-1,2){6}}
\put(30,12){\line(1,2){6}}
\put(36,24){\line(-1,1){12}}
\put(24,36){\line(2,1){24}}
\put(48,12){\line(-1,1){12}}
\put(42,0){\line(1,2){6}}
\put(36,-12){\line(1,2){6}}
\put(30,-24){\line(1,2){6}}
\put(36,-36){\line(-1,2){6}}
\put(0,0){\circle*{3}}
\put(6,12){\circle*{3}}
\put(18,12){\circle*{3}}
\put(12,24){\circle{3}}
\put(30,-12){\circle*{3}}
\put(24,0){\circle{3}}
\put(36,0){\circle*{3}}
\put(30,12){\circle{3}}
\put(36,24){\circle{3}}
\put(24,36){\circle{3}}
\put(30,-24){\circle{3}}
\put(48,12){\circle*{3}}
\put(42,0){\circle*{3}}
\put(36,-12){\circle*{3}}
\put(36,-36){\circle*{3}}
\put(48,48){\circle*{3}}
\put(48,48){\circle*{2}}
\put(48,12){\circle*{2} }
\put(42,0){\circle*{2}}
\put(40,-2){\dashbox(4,4){}} %\put(48,0){\makebox(0,0)[l]{$w=1221$}}
\put(40,-12){\vector(1,2){5}}
\qbezier[20](51,51)(51,30)(51,9)
\qbezier[20](51,9)(42,9)(38,1)
\qbezier[20](38,1)(36,-3)(30,-15)
\qbezier[20](30,-15)(-3,-15)(-3,-3)
\qbezier[40](-3,-3)(-3,51)(51,51)
\qbezier[20](29,-21)(32,-15)(34,-11)
\qbezier[6](34,-11)(35,-9)(36,-9)
\qbezier[6](36,-9)(39,-9)(39,-12)
\qbezier[40](39,-12)(39,-24)(39,-36)
\qbezier[6](39,-36)(39,-39)(36,-39)
\qbezier[6](36,-39)(35,-39)(34,-37)
\qbezier[20](34,-37)(32,-33)(29,-27)
\qbezier[10](29,-27)(27,-24)(29,-21)
 \end{picture}},
&T \nodeminus w =&
\raisebox{-84\unitlength}{\begin{picture}(58, 84)(96, -36)
\put(100,0){\line(1,2){6}}
\put(106,12){\line(1,2){6}}
\put(118,12){\line(-1,2){6}}
\put(112,24){\line(1,1){12}}
\put(130,-12){\line(-1,2){6}}
\put(124,0){\line(1,2){6}}
\put(136,0){\line(-1,2){6}}
\put(130,12){\line(1,2){6}}
\put(136,24){\line(-1,1){12}}
\put(124,36){\line(2,1){24}}
\put(148,12){\line(-1,1){12}}
\put(142,0){\line(1,2){6}}
\put(136,-12){\line(1,2){6}}
\put(142,-24){\line(-1,2){6}}
\put(100,0){\circle*{3}}
\put(106,12){\circle*{3}}
\put(118,12){\circle*{3}}
\put(112,24){\circle{3}}
\put(130,-12){\circle*{3}}
\put(124,0){\circle{3}}
\put(136,0){\circle*{3}}
\put(130,12){\circle{3}}
\put(136,24){\circle{3}}
\put(124,36){\circle{3}}
\put(136,-12){\circle{3}}
\put(148,12){\circle*{3}}
\put(142,0){\circle*{3}}
\put(142,-24){\circle*{3}}
\put(148,48){\circle*{3}}
\put(148,48){\circle*{2}}
\put(148,12){\circle*{2} }
\qbezier[20](151,51)(151,30)(151,9)
\qbezier[20](151,9)(142,9)(138,1)
\qbezier[20](138,1)(136,-3)(130,-15)
\qbezier[20](130,-15)(97,-15)(97,-3)
\qbezier[40](97,-3)(97,51)(151,51)
\qbezier[20](135,-9)(138,-3)(140,1)
\qbezier[6](140,1)(141,3)(142,3)
\qbezier[6](142,3)(145,3)(145,0)
\qbezier[40](145,0)(145,-12)(145,-24)
\qbezier[6](145,-24)(145,-27)(142,-27)
\qbezier[6](142,-27)(141,-27)(140,-25)
\qbezier[20](140,-25)(138,-21)(135,-15)
\qbezier[10](135,-15)(133,-12)(135,-9)
 \end{picture}},
\end{align*}
where $\bullet$ are nodes in $R_T$ or $R_{T\nodeminus w}$,
and 
\begin{picture}(4,4)
\put(0,0){\dashbox(4,4){\put(2,0){\circle*{3}}}}
\end{picture}
is $w=1221$.
The natural inclusion $\nu_{T,w}$ maps the nodes in 
\raisebox{0\unitlength}{\begin{picture}(14, 30)(26, -39)
\qbezier[20](29,-21)(32,-15)(34,-11)
\qbezier[6](34,-11)(35,-9)(36,-9)
\qbezier[6](36,-9)(39,-9)(39,-12)
\qbezier[40](39,-12)(39,-24)(39,-36)
\qbezier[6](39,-36)(39,-39)(36,-39)
\qbezier[6](36,-39)(35,-39)(34,-37)
\qbezier[20](34,-37)(32,-33)(29,-27)
\qbezier[10](29,-27)(27,-24)(29,-21)
 \end{picture}} of $T \nodeminus w$ to
the nodes in
\raisebox{0\unitlength}{\begin{picture}(14, 30)(26, -39)
\qbezier[20](29,-21)(32,-15)(34,-11)
\qbezier[6](34,-11)(35,-9)(36,-9)
\qbezier[6](36,-9)(39,-9)(39,-12)
\qbezier[40](39,-12)(39,-24)(39,-36)
\qbezier[6](39,-36)(39,-39)(36,-39)
\qbezier[6](36,-39)(35,-39)(34,-37)
\qbezier[20](34,-37)(32,-33)(29,-27)
\qbezier[10](29,-27)(27,-24)(29,-21)
 \end{picture}} of $T$, and the nodes in
\raisebox{0\unitlength}{\begin{picture}(56,66)(-4, -3)
\qbezier[20](51,51)(51,30)(51,9)
\qbezier[20](51,9)(42,9)(38,1)
\qbezier[20](38,1)(36,-3)(30,-15)
\qbezier[20](30,-15)(-3,-15)(-3,-3)
\qbezier[40](-3,-3)(-3,51)(51,51)
 \end{picture}} 
of $T \nodeminus w$
 to the nodes in
\raisebox{0\unitlength}{\begin{picture}(56,66)(-4, -3)
\qbezier[20](51,51)(51,30)(51,9)
\qbezier[20](51,9)(42,9)(38,1)
\qbezier[20](38,1)(36,-3)(30,-15)
\qbezier[20](30,-15)(-3,-15)(-3,-3)
\qbezier[40](-3,-3)(-3,51)(51,51)
 \end{picture}} of $T$, respectively.
}
\newcommand{\exampleiii}{
Let $T$ be as in Example \ref{exbintree}. Then $E_T$, $w_{T,k}$ and so on are the following:
\begin{gather*}
\begin{picture}(58, 84)(5, -36)
\put(0,0){\line(1,2){6}}
\put(6,12){\line(1,2){6}}
\put(18,12){\line(-1,2){6}}
\put(12,24){\line(1,1){12}}
\put(30,-12){\line(-1,2){6}}
\put(24,0){\line(1,2){6}}
\put(36,0){\line(-1,2){6}}
\put(30,12){\line(1,2){6}}
\thicklines
\put(36,24){\line(-1,1){12}}
\put(24,36){\line(2,1){24}}
\put(48,12){\line(-1,1){12}}
\put(42,0){\line(1,2){6}}
\put(36,-12){\line(1,2){6}}
\put(30,-24){\line(1,2){6}}
\put(36,-36){\line(-1,2){6}}
\thinlines
\put(0,0){\circle*{3}}
\put(6,12){\circle*{3}}
\put(18,12){\circle*{3}}
\put(12,24){\circle{3}}
\put(30,-12){\circle*{3}}
\put(24,0){\circle{3}}
\put(36,0){\circle*{3}}
\put(30,12){\circle{3}}
\put(36,24){\circle{3}}
\put(24,36){\circle{3}}
\put(30,-24){\circle{3}}
\put(48,12){\circle*{3}}
\put(42,0){\circle*{3}}
\put(36,-12){\circle*{3}}
\put(36,-36){\circle*{3}}
\put(48,48){\circle*{3}}
[6~%
\put(51,48){\makebox(0,0)[l]{$w_{T,1}$}}
\put(51,12){\makebox(0,0)[l]{$w_{T,2}$}}
\put(45,0){\makebox(0,0)[l]{$w_{T,3}$}}
\put(39,-12){\makebox(0,0)[l]{$w_{T,4}$}}
\put(39,-36){\makebox(0,0)[l]{$w_{T,5}$}}
 \end{picture}
,
\end{gather*}
where nodes on thick lines are in $E_T$ and $\bullet$ are in $R_T$.
}

\begin{document}

\title[An example of generalized Schur operators]{ An example of \\generalized Schur operators \\involving planar binary trees}
%\thanks{Department of Mathematics, Hokkaido University, Kita 10, Nishi 8, Kita-Ku, Sapporo,
%Hokkaido, 060-0810, Japan.  E-mail: {\tt nu@math.sci.hokudai.ac.jp}}
\author[Numata, y.]{Numata yasuhide}
\address{Kita 10, Nishi 8, Kita-Ku, Sapporo, Hokkaido, 060-0810, Japan.}
\email{nu@math.sci.hokudai.ac.jp}
%\urladdr{http://www.math.sci.hokudai.ac.jp/\textasciitilde nu/}
%\dedicatory
%\keywords{Robinson-Schensted algorithm; differential poset; 
%dual graded graph;
%Young tableaux.}
%\subjclass
\begin{abstract}
Young's lattice is a prototypical example of differential posets.
Differential posets have the Robinson correspondence, 
the correspondence between permutations and
pairs of standard tableaux with the same shape,
as in the case of Young's lattice. 
Fomin introduced
generalized Schur operators 
 to generalize 
the method of 
Robinson correspondence in differential posets
to the Robinson-Schensted-Knuth correspondence,
the correspondence between
 certain matrices and pairs of semi-standard tableaux with the same shape.
In this paper, we introduce operators on the 
vector space whose basis is the set of
planar binary trees.
To prove that the operators are generalized Schur operators,
we construct a correspondence, which is 
an extension of Fomin's $r$-correspondence for them.
\end{abstract}

\maketitle

\section{Introduction}
Stanley introduced differential posets in \cite{StaD, StaV}.
Young's lattice is a prototypical example of differential posets.
Young's lattice has  the Robinson correspondence, 
the correspondence between permutations and
pairs of standard tableaux whose shapes are the same Young diagram.
This correspondence was generalized for differential posets or
dual graphs (generalizations of differential posets
\cite{FomD}) 
by Fomin
\cite{FomRSK, FomS}. (See also \cite{Rob}.)

Young's lattice also has  the Robinson-Schensted-Knuth correspondence,
the correspondence between
 certain matrices and pairs of semi-standard tableaux.
Fomin \cite{FomK}
introduced operators called generalized Schur operators, and
generalized the method of the Robinson correspondence
to that of the Robinson-Schensted-Knuth correspondence.

In this paper, we introduce linear operators on 
the vector space whose basis is the  set of
binary trees.
We show the operators are generalized Schur operators.
To prove this, we construct an extension of
an $r$-correspondence.

\section{Generalized Schur Operators}\label{GSdefsec}

%First 
In this section,
we recall generalized Schur  operators introduced by Fomin \cite{FomK}.

Let $K$ be a field of characteristic zero
that contains all formal power series of variables
$t, t', t_1,t_2,\ldots$ %.
Let 
$V_i$ be finite-dimensional $K$-vector spaces for all $i \in \ZZ$. 
%and $V_i$ be zero for any $i<0$.
Fix a basis $Y_i$ of each $V_i$
so that $V_i=KY_i$.
%Let $Y=\bigcup_i Y_i$ and $V=KY$. 
Let $Y=\bigcup_i Y_i$,  $V=\bigoplus_i V_i$ and $\widehat{V}=\prod_{i} V_i$,  
i.e., $V$ is the vector space consisting of all finite linear
combinations of elements of $Y$, and
$\widehat{V}$ is the vector space consisting of all linear combinations
of elements of $Y$.
%
%Let $V$, $Y$ and $Y^\ast$ be $\bigoplus_i V_i$,
% $\bigcup_i Y_i$ 
% and the dual basis of $Y$, respectively.
%The \defit{rank function} on $V$  mapping $v \in V_i$ to $i$ 
%is denoted by $\rho$. 
%We say that $Y$  \defit{has a minimum} $\minimum$
%if $Y_i=\emptyset$ for $i<0$ and $Y_0=\{\minimum\}$.

For a sequence $\{A_i\}$ and a formal variable $x$, 
we write $A(x)$ for the generating function $\sum_{i\geq 0} A_i x^i$.

%Throughout this paper, for $i>0$,
%let $D_i$ and $U_i$ be linear operators on $V$ satisfying 
%$\rho(U_i\lambda)=\rho(\lambda)+i$ and 
%$\rho(D_i\lambda)=\rho(\lambda)-i$
%for $\lambda\in Y$.
%In other words, 
%the images $D_j(V_i)$ and $U_j(V_i)$ of $V_i$ by $D_j$ and $U_j$ are contained in $V_{i-j}$
%and $V_{i+j}$ for $i \in \ZZ$ and $j \in \N$ respectively.
%We call $D_i$ or $D(t)$ and $U_i$ or $U(t)$ \defit{down operators} and
%\defit{up operators}. 

\begin{definition}
We call $D(t_1)\cdots D(t_n)$ and $U(t_n)\cdots U(t_1)$
\defit{generalized Schur operators with} $\{a_m\}$
if the following conditions are satisfied$:$
\begin{itemize}
\item $\{a_m\}$ is a sequence of elements of $K$.
\item $U_i$ is a linear map on $V$ satisfying
$U_i(V_j) \subset V_{j+i}$
for all $j$.
\item $D_i$ is a linear map on $V$ satisfying
$D_i(V_j) \subset V_{j-i}$
for all $j$.
\item The equation 
$D(t')U(t)=a(t t')U(t)D(t')$
holds.
\end{itemize}
%
%Let $\{a_i\}$ be a sequence of elements of $K$.
%Down and up operators $D(t_1)\cdots D(t_n)$ and $U(t_n)\cdots U(t_1)$
%are called \defit{generalized Schur operators}
%if the equation $D(t')U(t)=a(t t')U(t)D(t')$
%holds. 
\end{definition}

%\begin{remark}
%In general, $D(t_1)\cdots D(t_n)$ and $U(t_n)\cdots U(t_1)$ are
%not linear operators on $V$  but  linear operators from $V$ to $\widehat{V}$.
%\end{remark}
\begin{remark}
Let  $\langle\phantom{x},\phantom{x}\rangle$ 
be the natural pairing  in $KY$,
i.e.,
the bilinear form 
on $\widehat{V} \times V$
such that 
$\langle\sum_{\lambda\in Y}a_\lambda \lambda,\sum_{\mu \in Y}b_\mu \mu\rangle =\sum_{\lambda \in Y}a_{\lambda} b_{\lambda}$.
For generalized Schur operators $D(t_1)\cdots D(t_n)$ and $U(t_n)\cdots U(t_1)$,
$U_i^{\ast}$ and $D_i^{\ast}$ denote
the maps obtained from the adjoints of $U_i$ and $D_i$ 
with respect to  $\langle\phantom{x},\phantom{x}\rangle$
by restricting to $V$, respectively.
For all $i$, 
$U_i^{\ast}$ and $D_i^{\ast}$ 
are linear maps on $V$
satisfying $U_i^{\ast}(V_j)\subset V_{j-i}$
and $D_i^{\ast}(V_{j})\subset V_{j+i}$.
By definition,
\begin{align*}
  \langle v, U_{i} w \rangle & = \langle w, U_{i}^{\ast} v \rangle,
&  \langle v, D_{i} w \rangle & = \langle w, D_{i}^{\ast} v \rangle
\end{align*}
for $v$, $w\in V$.
We write $U^{\ast}(t)$ and $D^{\ast}(t)$ for 
$\sum U^{\ast}_i t^i$ and
$\sum D^{\ast}_i t^i$.
By definition,
\begin{align*}
  \langle  U(t) \mu ,\lambda\rangle &= \langle  U^{\ast}(t) \lambda , \mu\rangle,
&  \langle  D(t) \mu ,\lambda\rangle &= \langle  D^{\ast}(t) \lambda , \mu\rangle
\end{align*}
for $\lambda$, $\mu\in Y$.
The equation  $D(t')U(t)=a(t t')U(t)D(t')$ implies the equation
$U^{\ast}(t')D^{\ast}(t)=a(t t')D^{\ast}(t)U^{\ast}(t')$.
Hence  $U^{\ast}(t_1)\cdots U^{\ast}(t_n)$ and 
$D^{\ast}(t_n)\cdots D^{\ast}(t_1)$ are 
generalized Schur operators with $\{a_m\}$ when 
 $D(t_1)\cdots D(t_n)$ and $U(t_n)\cdots U(t_1)$
are.% generalized Schur operators with $\{a_m\}$. 
\end{remark}

\section{Definition}
In this section,
 first we recall the definition of rooted planar binary trees and
labellings on them.
Next we introduce linear operators on
the vector space whose basis is the set of rooted planar binary trees.

\subsection{Rooted Planar Binary Trees}
We define rooted planar binary trees and their labellings.

Let $F$ be the monoid of words generated by the alphabet $\{1,2\}$
and let $0$ denote the word whose length is $0$.
We identify $F$ with a poset by $v \leq vw$ for $v,w \in F$.
We call an ideal of poset $F$ a \defit{rooted planar binary tree} or shortly \defit{tree}.
Let $\Trees$ denote the set of trees. 

Let $T$ be a tree.
An element of $T$ is called a \defit{node} of $T$.
We write $\Trees_i$ for the set of trees of $i$ nodes.
We respectively call nodes $v2$ and $v1$ \defit{right} and \defit{left
children} of $v$. 
A node without children is called a \defit{leaf}.
If $T$ is nonempty, $0\in T$. We call $0$ the \defit{root} of $T$.

For $T\in \Trees$ and $v \in F$,
we define $T_v$ by $T_v \denotes \Set{w\in T| v\leq w}$.

\begin{example}
For a tree $T=\Set{0,11,2,21,211,22,221,2211}$, the root, leaves
 and so on are as follows:
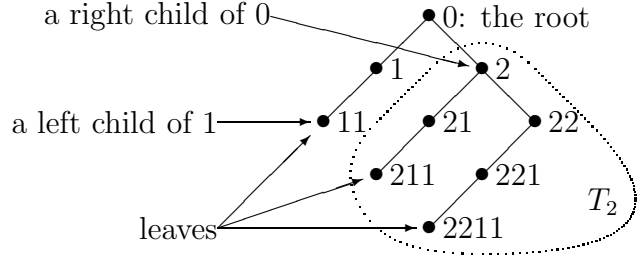
\begin{gather*}
\begin{picture}(200,100)(-100,-15)
\put(40,80){\line(1,-1){40}} 
\multiput(0,40)(20,-20){3}{\line(1,1){40}}
\multiput(0,40)(20,-20){3}{\makebox(0,0)[c]{$\bullet$}}
\multiput(20,60)(20,-20){3}{\makebox(0,0)[c]{$\bullet$}}
\multiput(40,80)(20,-20){3}{\makebox(0,0)[c]{$\bullet$}}
\put(45,80){\makebox(0,0)[l]{$0$: the root}}
\put(65,60){\makebox(0,0)[l]{$2$}}
\put(85,40){\makebox(0,0)[l]{$22$}}
\put(25,60){\makebox(0,0)[l]{$1$}}
\put(45,40){\makebox(0,0)[l]{$21$}}
\put(65,20){\makebox(0,0)[l]{$221$}}
\put(5,40){\makebox(0,0)[l]{$11$}}
\put(25,20){\makebox(0,0)[l]{$211$}}
\put(45,0){\makebox(0,0)[l]{$2211$}}
\qbezier[15](20,40)(30,50)(40,60)
\qbezier[15](80,60)(60,80)(40,60)
\qbezier[15](80,60)(90,50)(100,40)
\qbezier[60](60,-10)(150,-10)(100,40)
\qbezier[15](20,0)(30,-10)(60,-10)
\qbezier[15](20,0)(0,20)(20,40)
\put(100,10){\makebox(0,0)[l]{$T_2$}}
\put(-40,0){\vector(1,1){35}}
\put(-40,0){\vector(3,1){55}}
\put(-40,0){\vector(1,0){75}}
\put(-40,0){\makebox(0,0)[r]{leaves}}
\put(-40,40){\vector(1,0){35}}
\put(-40,40){\makebox(0,0)[r]{a left child of $1$}}
\put(-20,80){\vector(4,-1){75}}
\put(-20,80){\makebox(0,0)[r]{a right child of $0$}}
\end{picture}.
\end{gather*}

\end{example}

\begin{definition}
Let $T$ be a tree and $m$ a positive integer.
We call a map $\varphi: T \to \Set{1,\ldots , m}$
a right-strictly-increasing labelling if 
\begin{itemize}
\item $\varphi(w)\leq \varphi(v)$ for $w \in T$ and $v \in T_{w1}$ and
\item $\varphi(w)< \varphi(v)$ for $w \in T$ and $v \in T_{w2}$.
\end{itemize}
We call a map $\phi: T \to \Set{1,\ldots , m}$
a left-strictly-increasing labelling if 
\begin{itemize}
\item $\phi(w)< \phi(v)$ for $w \in T$ and $v \in T_{w1}$ and
\item $\phi(w)\leq \phi(v)$ for $w \in T$ and $v \in T_{w2}$.
\end{itemize}
We call a map $\psi: T \to \Set{1,\ldots , m}$
a binary-searching labelling if 
\begin{itemize}
\item $\psi(w)\geq \psi(v)$ for $w \in T$ and $v \in T_{w1}$ and
\item $\psi(w)< \psi(v)$ for $w \in T$ and $v \in T_{w2}$.
\end{itemize}
\end{definition}

\begin{example}
The following is a right-strictly-increasing labelling:
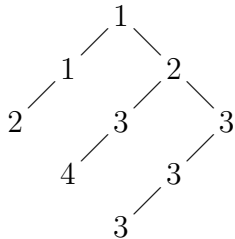
\begin{gather*}
\begin{picture}(100,100)
\multiput(5,45)(20,-20){3}{\line(1,1){10}}
\multiput(25,65)(20,-20){3}{\line(1,1){10}}
\multiput(45,75)(20,-20){2}{\line(1,-1){10}}
\put(40,80){\makebox(0,0)[c]{$1$}}
\put(60,60){\makebox(0,0)[c]{$2$}}
\put(80,40){\makebox(0,0)[c]{$3$}}
\put(20,60){\makebox(0,0)[c]{$1$}}
\put(40,40){\makebox(0,0)[c]{$3$}}
\put(60,20){\makebox(0,0)[c]{$3$}}
\put(0,40){\makebox(0,0)[c]{$2$}}
\put(20,20){\makebox(0,0)[c]{$4$}}
\put(40,0){\makebox(0,0)[c]{$3$}}
\end{picture}.
\end{gather*}
The following is a left-strictly-increasing labelling:
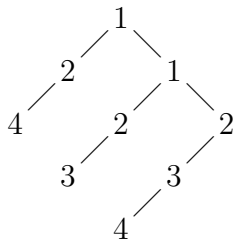
\begin{gather*}
\begin{picture}(100,100)
\multiput(5,45)(20,-20){3}{\line(1,1){10}}
\multiput(25,65)(20,-20){3}{\line(1,1){10}}
\multiput(45,75)(20,-20){2}{\line(1,-1){10}}
\put(40,80){\makebox(0,0)[c]{$1$}}
\put(60,60){\makebox(0,0)[c]{$1$}}
\put(80,40){\makebox(0,0)[c]{$2$}}
\put(20,60){\makebox(0,0)[c]{$2$}}
\put(40,40){\makebox(0,0)[c]{$2$}}
\put(60,20){\makebox(0,0)[c]{$3$}}
\put(0,40){\makebox(0,0)[c]{$4$}}
\put(20,20){\makebox(0,0)[c]{$3$}}
\put(40,0){\makebox(0,0)[c]{$4$}}
\end{picture}.
\end{gather*}
The following is a binary-searching labelling:
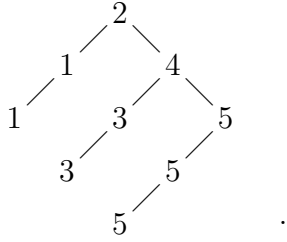
\begin{gather*}
\begin{picture}(100,100)
\multiput(5,45)(20,-20){3}{\line(1,1){10}}
\multiput(25,65)(20,-20){3}{\line(1,1){10}}
\multiput(45,75)(20,-20){2}{\line(1,-1){10}}
\put(40,80){\makebox(0,0)[c]{$2$}}
\put(60,60){\makebox(0,0)[c]{$4$}}
\put(80,40){\makebox(0,0)[c]{$5$}}
\put(20,60){\makebox(0,0)[c]{$1$}}
\put(40,40){\makebox(0,0)[c]{$3$}}
\put(60,20){\makebox(0,0)[c]{$5$}}
\put(0,40){\makebox(0,0)[c]{$1$}}
\put(20,20){\makebox(0,0)[c]{$3$}}
\put(40,0){\makebox(0,0)[c]{$5$}}
\end{picture}.
\end{gather*}

\end{example}

\subsection{Definition of our generalized Schur operators}
In this section, we define linear operators $U_i$, $U'_i$, $D_i$.
In Section \ref{mainsec},
we shall show that these linear operators are generalized Schur operators.

\subsubsection{Up operators}
First we define up operators $U_i$ 
and consider a relation between the up operators $U_i$
and 
right-strictly labellings.
Next we define $U'_i$ and consider a relation between 
the up operators $U'_i$ and 
left-strictly labellings.

\begin{definition}
We define the edges $G_{U}$ of oriented graphs whose
 vertices are trees to be the set of pairs $(T,T')$ of trees satisfying the
 following:
\begin{itemize}
\item $T\subset T'$.
\item For each $w \in T' \setminus T$, there exists $v_w \in T$ such that
      $w=v_w1^n$ or  $w=v_w21^n$ for some nonnegative integer $n$ if
      $T \neq \emptyset$.
\item For each $w \in T' \setminus T$, $w=1^n$ for some nonnegative
      integer $n$ if $T= \emptyset$.
\end{itemize}
We call $T'$ \defit{a tree obtained from $T$ by adding some nodes right-strictly}
if $(T,T')\in G_{U}$.
We define $G_{U_i}$ by
\begin{gather*}
G_{U_i}=\Set{(T,T') \in G_{U} | \numof{T}+i= \numof{T'}}.
\end{gather*}
\end{definition}

\begin{definition}\label{defofu}
For $i\in \NN$ and $T\in \Trees$, we define  linear operators $U_i$ on $K \Trees$ by
\begin{gather*}
U_i T = \sum_{T':(T,T')\in G_{U_i}} T' .
\end{gather*}
Equivalently, $U_i T$ is the sum of all trees obtained from $T$ by adding $i$ nodes
right-strictly. 
\end{definition}

\begin{example}
For example, $U_3$ acts on $\Set{0}$ as follows:
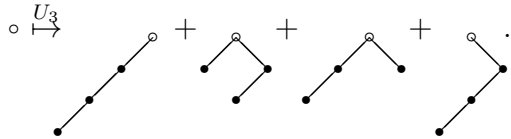
\begin{gather*}
\begin{picture}(6, 6)(-3,-3)
\put(0,0){\circle{3}}
\end{picture}
\overset{U_3}{\mapsto}
\raisebox{-36\unitlength}{\begin{picture}(36, 36)(5,0)
\put(0,0){\line(1,1){12}}
\put(12,12){\line(1,1){12}}
\put(24,24){\line(1,1){12}}
\put(36,36){\circle{3}}
%\color{\addcolor}
\put(24,24){\circle*{3}}
\put(0,0){\circle*{3}}
\put(12,12){\circle*{3}}
%\color{black}
\end{picture}}
+
\raisebox{-24\unitlength}{\begin{picture}(24, 24)(0,-12)
\put(0,0){\line(1,1){12}}
\put(12,-12){\line(1,1){12}}
\put(24,0){\line(-1,1){12}}
\put(12,12){\circle{3}}
%\color{\addcolor}
\put(0,0){\circle*{3}}
\put(24,0){\circle*{3}}
\put(12,-12){\circle*{3}}
%\color{black}
\end{picture}}
+
\raisebox{-24\unitlength}{\begin{picture}(36, 24)(-12,-12)
\put(-12,-12){\line(1,1){12}}
\put(0,0){\line(1,1){12}}
\put(24,0){\line(-1,1){12}}
\put(12,12){\circle{3}}
%\color{\addcolor}
\put(0,0){\circle*{3}}
\put(24,0){\circle*{3}}
\put(-12,-12){\circle*{3}}
%\color{black}
\end{picture}}
+
\raisebox{-36\unitlength}{\begin{picture}(24, 36)(-12,-12)
\put(-12,-12){\line(1,1){12}}
\put(12,12){\line(-1,1){12}}
\put(0,0){\line(1,1){12}}
\put(0,24){\circle{3}}
%\color{\addcolor}
\put(12,12){\circle*{3}}
\put(0,0){\circle*{3}}
\put(-12,-12){\circle*{3}}
%\color{black}
\end{picture}}.
\end{gather*}
\end{example}

\begin{remark}
Let $\varphi$ be a right-strictly-increasing labelling.
The inverse image $\varphi^{-1}(\{1,\ldots,n+1\})$ 
is the tree obtained
from the inverse image $\varphi^{-1}(\{1,\ldots,n\})$
by adding  some nodes right-strictly.
Hence we identify right-strictly-increasing labellings 
with paths $(\emptyset=T^{0},T^{1},\ldots, T^{m})$ of $G_{U}$.
\end{remark}

\begin{example}
For example, we identify a right-strictly-increasing labelling
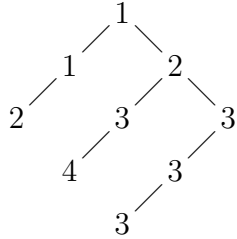
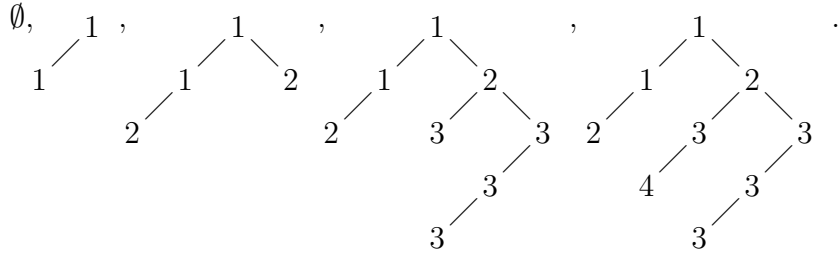
\begin{gather*}
\begin{picture}(100,100)
\multiput(5,45)(20,-20){3}{\line(1,1){10}}
\multiput(25,65)(20,-20){3}{\line(1,1){10}}
\multiput(45,75)(20,-20){2}{\line(1,-1){10}}
\put(40,80){\makebox(0,0)[c]{$1$}}
\put(60,60){\makebox(0,0)[c]{$2$}}
\put(80,40){\makebox(0,0)[c]{$3$}}
\put(20,60){\makebox(0,0)[c]{$1$}}
\put(40,40){\makebox(0,0)[c]{$3$}}
\put(60,20){\makebox(0,0)[c]{$3$}}
\put(0,40){\makebox(0,0)[c]{$2$}}
\put(20,20){\makebox(0,0)[c]{$4$}}
\put(40,0){\makebox(0,0)[c]{$3$}}
\end{picture}
\end{gather*}
with a sequence
\begin{gather*}
\emptyset,
\raisebox{-80\unitlength}{\begin{picture}(30,80)(20,0)
%\multiput(5,45)(20,-20){3}{\line(1,1){10}}
\multiput(25,65)(20,-20){1}{\line(1,1){10}}
%\multiput(45,75)(20,-20){2}{\line(1,-1){10}}
%
\put(40,80){\makebox(0,0)[c]{$1$}}
\put(20,60){\makebox(0,0)[c]{$1$}}
\end{picture}}
,
\raisebox{-80\unitlength}{\begin{picture}(70,80)
\multiput(5,45)(20,-20){1}{\line(1,1){10}}
\multiput(25,65)(20,-20){1}{\line(1,1){10}}
\multiput(45,75)(20,-20){1}{\line(1,-1){10}}
\put(40,80){\makebox(0,0)[c]{$1$}}
\put(60,60){\makebox(0,0)[c]{$2$}}
\put(20,60){\makebox(0,0)[c]{$1$}}
\put(0,40){\makebox(0,0)[c]{$2$}}
\end{picture}}
,
\raisebox{-80\unitlength}{\begin{picture}(90,80)
\multiput(5,45)(40,-40){2}{\line(1,1){10}}
\multiput(25,65)(20,-20){3}{\line(1,1){10}}
\multiput(45,75)(20,-20){2}{\line(1,-1){10}}
\put(40,80){\makebox(0,0)[c]{$1$}}
\put(60,60){\makebox(0,0)[c]{$2$}}
\put(80,40){\makebox(0,0)[c]{$3$}}
\put(20,60){\makebox(0,0)[c]{$1$}}
\put(40,40){\makebox(0,0)[c]{$3$}}
\put(60,20){\makebox(0,0)[c]{$3$}}
\put(0,40){\makebox(0,0)[c]{$2$}}
\put(40,0){\makebox(0,0)[c]{$3$}}
\end{picture}}
,\ 
\raisebox{-80\unitlength}{\begin{picture}(90,80)
\multiput(5,45)(20,-20){3}{\line(1,1){10}}
\multiput(25,65)(20,-20){3}{\line(1,1){10}}
\multiput(45,75)(20,-20){2}{\line(1,-1){10}}
\put(40,80){\makebox(0,0)[c]{$1$}}
\put(60,60){\makebox(0,0)[c]{$2$}}
\put(80,40){\makebox(0,0)[c]{$3$}}
\put(20,60){\makebox(0,0)[c]{$1$}}
\put(40,40){\makebox(0,0)[c]{$3$}}
\put(60,20){\makebox(0,0)[c]{$3$}}
\put(0,40){\makebox(0,0)[c]{$2$}}
\put(20,20){\makebox(0,0)[c]{$4$}}
\put(40,0){\makebox(0,0)[c]{$3$}}
\end{picture}}.
\end{gather*}

\end{example}

Next we define other up operators $U'_{i}$.
\begin{definition}
We define the edges $G_{U'}$  of oriented graphs whose
 vertices are trees to be the set of pairs $(T,T')$ of trees satisfying the
 following: 
\begin{itemize}
\item $T\subset T'$.
\item For each $w \in T' \setminus T$, there exists $v_w \in T$ such that
      $w=v_w2^n$ or  $w=v_w12^n$ for some nonnegative integer $n$ if
      $T\neq\emptyset$.
\item For each $w \in T' \setminus T$, 
      $w=2^n$ for some nonnegative integer $n$ if
      $T=\emptyset$.
\end{itemize}
We call $T'$ \defit{a tree obtained from $T$ by adding some nodes
 left-strictly} if $(T,T')\in G_{U'}$.
We define $G_{U'_i}$ by
\begin{gather*}
G_{U'_i} = \Set{(T,T')\in G_{U}|\numof{T}+i=\numof{T'}}.
\end{gather*} 
\end{definition}
\begin{definition}\label{defofdualu}
For $i\in \NN$ and $T\in \Trees$, we define  linear operators $U'_i$ on $K \Trees$ to be
\begin{gather*}
U'_i T = \sum_{T': (T,T')\in G_{U'_i}} T' .
\end{gather*}
Equivalently, $U'_i T$ is the sum of all trees obtained from $T$ by adding $i$ nodes
left-strictly. 
\end{definition}

\begin{remark}
Similarly to the case of $U_i$ and right-strictly-increasing labellings, 
we identify left-strictly-increasing labellings  
with paths $(\emptyset=T^{0},T^{1},\ldots, T^{m})$ of $G_{U'}$.
\end{remark}

\subsubsection{Down operators}
Next we define down operators $D_i$  on $K\Trees$ and we see relations
between the down operators $D_i$ and binary searching labellings.

First we prepare some terms to define the linear operators $D_i$.

For $T \in \Trees$,
let $R_T$ denote the set $\{w \in T | w2 \not\in T\}$,
i.e., the set of nodes of $T$ without right children.
For $w \in R_T$, we define
\begin{gather*}
T \nodeminus w = (T \setminus T_w) \cup \set{wv|w1v \in T_w}.
\end{gather*}
There exists the natural inclusion $\nu_{T , w}$ from $T\nodeminus w$
to $T$ defined by
\begin{gather*}
\begin{cases}
\nu_{T , w}(wv)=w1v & (wv \in T_{w})\\
\nu_{T , w}(v')=v' & (v' \not\in T_{w}). 
\end{cases}
\end{gather*}
\begin{example}\label{exbintree}
\exampleii
\end{example}
For $T\in \Trees$, let $E_T$ denote $\{ w \in T |$ 
If $w=v1w'$ then $v2 \not\in T \}$.
Roughly speaking,
it is the set of nodes of $T$ between the root $0$ and the right-most node of $T$.
We define  $r_T$ by  $r_T = E_T \cap R_T$. 
The set $r_T$ is a chain.
Let $r_T=\{w_{T,1} < w_{T,2} < w_{T,3} < \cdots < w_{T,k}\}$.
Let 
$r_{T,i}$ 
denote 
the ideal
$\{w_{T,1}, w_{T,2}, w_{T,3}, \ldots, w_{T,i}\}$ 
of $r_T$
consisting of $i$ nodes.
\begin{example}
\exampleiii
\end{example}

We define $T \stripminus r_{T,i}$ inductively by
\begin{gather*}
%(\ldots((T \nodeminus w_{T,i}) \nodeminus w_{T,i-1})\nodeminus \cdots )\nodeminus r_{T,1}.
\begin{cases}
(T\nodeminus w_{T,i}) \stripminus r_{T,i-1}  & i > 0\\
T & i=0.
\end{cases}
\end{gather*}
The natural inclusion $\nu_{T , w}$ 
induces the natural inclusion 
\begin{gather*}
%\nu_{T,i}=\nu_{T,i}\circ\nu_{T,i}\circ \cdots \circ 
%\nu_{T\stripminus r_{T,i-1},w{T,i-1}} 
%\nu_{T,r_{T,i}}=  \nu_{T\nodeminus w_{T,i}, r_{T,i-1}} \circ \nu_{T, w_{T,i}}
\nu_{T,i}=  \nu_{T \nodeminus w_{T,i},i-1} \circ \nu_{T, w_{T,i}}
\end{gather*}
from $T\stripminus r_{T,i}$
to $T$. 
\begin{example}
\examplei
\end{example}

We also define a bijection $\widetilde \nu_{T,i}$ from the words $F$ of $\{1,2\}$ to
$F\setminus r_{T,i}$
\begin{gather*}
\widetilde \nu_{T,i}(w) = \nu_{T,i}(v)v',
\end{gather*}
where $w=vv'$ and $v=\max\set{u\in T\stripminus r_{T,i}|w=uu'}$.
By definition, $\widetilde \nu_{T,i}(w)= \nu_{T,i}(w)$
for $w\in T\stripminus r_{T,i}$.

\begin{definition}
We define the edges $G_{D_i}$ and $G_{D}$ of graphs whose vertices are
trees by
\begin{gather*}
G_{D_i}=\Set{ (T \stripminus r_{T,i} , T) | \numof{r_{T}} \geq i}
\intertext{and}
G_{D}=\bigcup_i G_{D_i}.
\end{gather*}
\end{definition}

\begin{remark}
By definition, $G_{D_0}=\Set{(T,T)|T\in \Trees}$.
For each $i$ and each $T\in \Trees$, the in-degree of $T$ in $G_{D_i}$, 
i.e. $\numof{\Set{(T',T)\in G_{D_i}}}$, 
is $1$.
\end{remark}

\begin{definition}\label{defofd}
We define  linear operators $D_i T$ to be 
$T'$ such that $(T',T)\in G_{D_i}$
%\begin{gather*}
%%\begin{cases}
%\sum_{T':(T',T)\in G_{D_i}} T'
%%T \stripminus r_{T,i} & i \leq \numof{r_T} \\
%%0 & i > \numof{r_T}
%%\end{cases}
%\end{gather*}
for $T\in\Trees$.
\end{definition}
Roughly speaking, $D_i T $ 
is the tree obtained from $T$ by 
evacuating the $i$ topmost nodes 
without a child on its right 
%and belonging 
between the root $0$ and the rightmost leaf of $T$.

\begin{example}
For example, $D_3$ acts as follows:
\begin{gather*}
\raisebox{-84\unitlength}{\begin{picture}(58, 84)(5, -36)
\put(0,0){\line(1,2){6}}
\put(6,12){\line(1,2){6}}
\put(18,12){\line(-1,2){6}}
\put(12,24){\line(1,1){12}}
\put(30,-12){\line(-1,2){6}}
\put(24,0){\line(1,2){6}}
\put(36,0){\line(-1,2){6}}
\put(30,12){\line(1,2){6}}
%
%\color{\lmlcolor}
%
\put(36,24){\line(-1,1){12}}
\put(24,36){\line(2,1){24}}
\put(48,12){\line(-1,1){12}}
\put(42,0){\line(1,2){6}}
\put(36,-12){\line(1,2){6}}
\put(30,-24){\line(1,2){6}}
\put(36,-36){\line(-1,2){6}}
%\color{black}
%%%%%%%%
\put(0,0){\circle{3}}
\put(6,12){\circle{3}}
\put(18,12){\circle{3}}
\put(12,24){\circle{3}}
\put(30,-12){\circle{3}}
\put(24,0){\circle{3}}
\put(36,0){\circle{3}}
\put(30,12){\circle{3}}
\put(36,24){\circle{3}}
\put(24,36){\circle{3}}
\put(30,-24){\circle{3}}
%
%\color{\worccolor}
\put(48,12){\circle{3}}
\put(42,0){\circle{3}}
\put(36,-12){\circle{3}}
\put(36,-36){\circle{3}}
\put(48,48){\circle{3}}
%%%%%%%%
%\color{\tmncolor}
\put(48,48){\circle*{2}}
\put(48,12){\circle*{2}}
\put(42,0){\circle*{2}}
%%%%%%
%\color{\evcolor}
\put(30,36){\vector(2,1){16}}
\put(40,-12){\vector(1,2){10}}
%
%\color{black}
 \end{picture}}
\overset{D_3}{\mapsto}%
\raisebox{-48\unitlength}{\begin{picture}(48, 60)(0, -12)
\put(0,0){\line(1,2){6}}
\put(6,12){\line(1,2){6}}
\put(18,12){\line(-1,2){6}}
\put(12,24){\line(1,1){12}}
\put(30,-12){\line(-1,2){6}}
\put(24,0){\line(1,2){6}}
\put(36,0){\line(-1,2){6}}
\put(30,12){\line(1,2){6}}
%
%\color{\lmlcolor}
\put(36,24){\line(-1,1){12}}
%
%\put(24,36){\line(2,1){24}}
%
\put(48,12){\line(-1,1){12}}
%
%\put(42,0){\line(1,2){6}}
%\put(36,-12){\line(1,2){6}}
\put(42,0){\line(1,2){6}}
\put(48,-12){\line(-1,2){6}}
%\color{black}
%%%%%%%%
\put(0,0){\circle{3}}
\put(6,12){\circle{3}}
\put(18,12){\circle{3}}
\put(12,24){\circle{3}}
\put(30,-12){\circle{3}}
\put(24,0){\circle{3}}
\put(36,0){\circle{3}}
\put(30,12){\circle{3}}
\put(36,24){\circle{3}}
\put(24,36){\circle{3}}
%
%\put(48,12){\circle{3}}
%
%\put(42,0){\circle{3}}
\put(42,0){\circle{3}}
%
%\color{\worccolor}
\put(48,12){\circle{3}}
\put(48,-12){\circle{3}}
%\put(48,48){\circle{3}}
%%%%%%%%
%\put(48,48){\circle*{2}}
%\put(48,12){\circle*{2}}
%\put(42,0){\circle*{2}}
 \end{picture}}.
\end{gather*}
\end{example}

Next we consider a relation between $G_{D}$ and
binary-searching labellings.

\begin{remark} 
Let $\psi_{m}:T \to \{1,\ldots, m\}$ be  a binary-searching labelling.
By the definition of binary-searching labelling,
the inverse image $\psi_{m}^{-1}(\{m\})$
equals $r_{T,k_m}=\{w_{T,1},\ldots,w_{T,k_m}\}$ for some $k_m$.
Hence we can construct $T \stripminus \psi_{m}^{-1}(\{m\})$.
The natural inclusion $\nu_{T, k_m}$
induces  a binary-searching labelling 
\begin{gather*}
\psi_{m-1}=\psi_{m}\circ \nu_{T, k_m}: 
T\stripminus \psi_{m}^{-1}(\{m\}) \to \{1,\ldots, m-1\}.
\end{gather*}
Hence we identify binary-searching labellings
with paths 
\begin{gather*}
(\emptyset=T^{0},T^{1},\ldots, T^{m})
\end{gather*}
 of $G_{D}$.
\end{remark}

\begin{example}
For example, we identify a binary-searching labelling
\begin{gather*}
\begin{picture}(100,100)
\multiput(5,45)(20,-20){3}{\line(1,1){10}}
\multiput(25,65)(20,-20){3}{\line(1,1){10}}
\multiput(45,75)(20,-20){2}{\line(1,-1){10}}
\put(40,80){\makebox(0,0)[c]{$2$}}
\put(60,60){\makebox(0,0)[c]{$4$}}
\put(80,40){\makebox(0,0)[c]{$5$}}
\put(20,60){\makebox(0,0)[c]{$1$}}
\put(40,40){\makebox(0,0)[c]{$3$}}
\put(60,20){\makebox(0,0)[c]{$5$}}
\put(0,40){\makebox(0,0)[c]{$1$}}
\put(20,20){\makebox(0,0)[c]{$3$}}
\put(40,0){\makebox(0,0)[c]{$5$}}
\end{picture}.
\end{gather*}
with a sequence
\begin{gather*}
\emptyset,
\raisebox{-80\unitlength}{\begin{picture}(30,80)(20,0)
%\multiput(5,45)(20,-20){3}{\line(1,1){10}}
\multiput(25,65)(20,-20){1}{\line(1,1){10}}
%\multiput(45,75)(20,-20){2}{\line(1,-1){10}}
%
\put(40,80){\makebox(0,0)[c]{$1$}}
\put(20,60){\makebox(0,0)[c]{$1$}}
\end{picture}},
\raisebox{-80\unitlength}{\begin{picture}(50,80)
\multiput(5,45)(20,-20){1}{\line(1,1){10}}
\multiput(25,65)(20,-20){1}{\line(1,1){10}}
%\multiput(45,75)(20,-20){2}{\line(1,-1){10}}
%
\put(40,80){\makebox(0,0)[c]{$2$}}
\put(20,60){\makebox(0,0)[c]{$1$}}
\put(0,40){\makebox(0,0)[c]{$1$}}
\end{picture}},
\raisebox{-80\unitlength}{\begin{picture}(70,80)
\multiput(5,45)(20,-20){1}{\line(1,1){10}}
\multiput(25,65)(20,-20){2}{\line(1,1){10}}
\multiput(45,75)(20,-20){1}{\line(1,-1){10}}
\put(40,80){\makebox(0,0)[c]{$2$}}
\put(60,60){\makebox(0,0)[c]{$3$}}
\put(20,60){\makebox(0,0)[c]{$1$}}
\put(40,40){\makebox(0,0)[c]{$3$}}
\put(0,40){\makebox(0,0)[c]{$1$}}
\end{picture}},
\raisebox{-80\unitlength}{\begin{picture}(70,80)
\multiput(5,45)(20,-20){2}{\line(1,1){10}}
\multiput(25,65)(20,-20){2}{\line(1,1){10}}
\multiput(45,75)(20,-20){1}{\line(1,-1){10}}
\put(40,80){\makebox(0,0)[c]{$2$}}
\put(60,60){\makebox(0,0)[c]{$4$}}
\put(20,60){\makebox(0,0)[c]{$1$}}
\put(40,40){\makebox(0,0)[c]{$3$}}
\put(0,40){\makebox(0,0)[c]{$1$}}
\put(20,20){\makebox(0,0)[c]{$3$}}
\end{picture}},
\raisebox{-80\unitlength}{\begin{picture}(80,80)
\multiput(5,45)(20,-20){3}{\line(1,1){10}}
\multiput(25,65)(20,-20){3}{\line(1,1){10}}
\multiput(45,75)(20,-20){2}{\line(1,-1){10}}
\put(40,80){\makebox(0,0)[c]{$2$}}
\put(60,60){\makebox(0,0)[c]{$4$}}
\put(80,40){\makebox(0,0)[c]{$5$}}
\put(20,60){\makebox(0,0)[c]{$1$}}
\put(40,40){\makebox(0,0)[c]{$3$}}
\put(60,20){\makebox(0,0)[c]{$5$}}
\put(0,40){\makebox(0,0)[c]{$1$}}
\put(20,20){\makebox(0,0)[c]{$3$}}
\put(40,0){\makebox(0,0)[c]{$5$}}
\end{picture}}.
\end{gather*}
\end{example}

\section{Main Results}\label{mainsec}

In this section, we show 
that  %operators  
 $D(t_1)\cdots D(t_n)$ and  $U(t_n)\cdots U(t_1)$
are generalized  Schur operators with $\{1,1,1,\ldots\}$.
We also show that 
%linear operators  
$D(t_1)\cdots D(t_n)$ and  
$U'(t_n)\cdots U'(t_1)$
are generalized  Schur operators with $\{1,1,0,0,0,0,\ldots\}$.
To prove the assertion for  $D(t_1)\cdots D(t_n)$ and  $U(t_n)\cdots U(t_1)$, 
we construct 
correspondences between  $N_{i,j}(T,T')$ and $\widetilde S_{j,i}(T,T')$, 
where $N_{i,j}(T,T')$ is the set of paths of graphs from $T$
to $T'$ via $G_{D_j}$ after $G_{U_i}$,
and $\widetilde S_{j,i}(T,T')$ is the set of paths of graphs from $T$
to $T'$ via $G_{U_{i-k}}$ after $G_{D_{j-k}}$,
where $o\leq k \leq \min\{i,j\}$. 
(See Definition \ref{defofrcorres}.)
The correspondence for $i=j=1$ is an $r$-correspondence introduced by
Fomin \cite{FomS},
which is needed to construct Robinson correspondences for $r$-dual graphs.
We also prove the assertion for  $D(t_1)\cdots D(t_n)$ and  $U(t_n)\cdots U(t_1)$ 
similarly.

\subsection{Main Theorems}
We prove the following theorems in Section \ref{proofsec}.

\begin{theorem}\label{mainthm}
Let $D_i$ be the linear operators defined in Definition \ref{defofd}
and $U_i$    the linear operators defined in Definition \ref{defofu}.
The operators  $D(t)$ and  $U(t')$
satisfy the equation
\begin{gather*}
D(t) U(t') = \frac{1}{1-t t'} U(t') D(t)  .
\end{gather*}
Equivalently, operators  $D(t_1)\cdots D(t_n)$ and  $U(t_n)\cdots U(t_1)$
are generalized  Schur operators with $\{1,1,1,\ldots\}$.
\end{theorem}

\begin{theorem}\label{dualmainthm}
Let $D_i$ be the linear operators defined in Definition \ref{defofd}
and $U'_i$    the linear operators defined in Definition \ref{defofdualu}.
The operators  $D(t)$ and  $U'(t')$
satisfy the equation
\begin{gather*}
D(t) U'(t') = (1 + t t') U'(t') D(t)  .
\end{gather*}
Equivalently, operators  $D(t_1)\cdots D(t_n)$ and  $U'(t_n)\cdots U'(t_1)$
are generalized  Schur operators with $\{1,1,0,0,0,0\ldots\}$.
\end{theorem}

\begin{cor}
A pair $(G_{U_1}=G_{U'_1}, G_{D_i})$
of graded graphs is an example of $1$-dual graphs %,
%where $1$-dual graphs means as in \cite{FomK}.
in the sense of Fomin \cite{FomK}.
Equivalently, $U_1$ and $D_1$ satisfy the equation
\begin{gather*}
D_1 U_1 - U_1 D_1 = I,
\end{gather*}
where $I$ is the identity map on $V$.
\end{cor}

\begin{remark}
Janvier Nzeutchap \cite{Nzc}
constructs $r$-dual graphs 
from dual Hopf  algebras.
The $1$-dual graphs obtained from the Loday-Ronco algebra
by his method are $G_{U_1}$ and $G_{D_1}$.
\end{remark}

\begin{cor}
The up and down operators $(U_1, D)$ of $(G_{U_1}, G_{D})$
satisfy
\begin{gather*}
D U_1  - U_1 D = D .
\end{gather*} 
\end{cor}

\begin{cor}
The up and down operators $( D_1 ^{\ast}, U^{\ast})$ 
of $(G_{D_1}, G_{U})$
satisfy
\begin{gather*}
U^{\ast} D_1^{\ast}  - D_1^{\ast} U^{\ast} = U^{\ast} .
\end{gather*} 
\end{cor}

\subsection{Proof of Main results}\label{proofsec} 
In this section, we prove Theorem \ref{mainthm}, i.e., 
\begin{gather*}
D(t) U(t') = \frac{1}{1-t t'} U(t') D(t),
\end{gather*}
and 
Theorem \ref{dualmainthm}, i.e., 
\begin{gather*}
D(t) U'(t') = (1 + t t') U'(t') D(t)  .
\end{gather*}
First we rewrite these equations as the equations of
the numbers of elements of some sets (Remark \ref{rewriterem}).
We show the equation by constructing  bijections (Lemmas
\ref{rcorreslemma} and \ref{dualrcorreslemma}).

\begin{lemma}
The equation
\begin{gather*}
D(t) U(t') = \frac{1}{1-t t'} U(t') D(t)
\end{gather*}
is equivalent to the equations
\begin{align}
D_j U_i= \sum_{k=0}^{\min{(i,k)}} U_{i-k} D_{j-k}\label{opeq1} 
& & \text{for all $i,j$.
}
\end{align}
\end{lemma}

\begin{lemma}
The equation
\begin{gather*}
D(t) U'(t') = {(1+t t')} U'(t') D(t)
\end{gather*}
is equivalent to the equations
\begin{align}
D_j U'_i&=  
%\begin{cases}
%U'_{i} D_{j} +  U'_{i-1} D_{j-1} & (i,j>0)\\
%U'_{i} D_{j} & \text{otherwise}
%\end{cases}
%&=
\sum_{k=0}^{\min{(1,i,k)}} U'_{i-k} D_{j-k}
\label{opeq2}
& \text{for all $i,j$.}
\end{align}
\end{lemma}

\begin{definition}\label{defofrcorres}
We respectively define sets $N_{i,j}(T,T')$ and $N'_{i,j}(T,T')$
of paths to be
\begin{gather*}
\Set{((T,T''),(T',T''))\in G_{U_i}\times G_{D_j}} \\
\intertext{and}
\Set{((T,T''),(T',T''))\in G_{U'_i}\times G_{D_j}}. \\
\end{gather*}
We respectively define sets $S_{j,i}(T,T')$ and  $S'_{j,i}(T,T')$
of paths to be
\begin{gather*}
\Set{((T'',T),(T'',T'))\in G_{D_j}\times G_{U_i}} \\
\intertext{and}
\Set{((T'',T),(T'',T'))\in G_{D_j}\times G_{U'_i}}.
\end{gather*}
We also define 
$\widetilde S_{j,i}(T,T')$ and  $\widetilde  S'_{j,i}(T,T')$
by 
\begin{gather*}
\widetilde S_{j,i}(T,T')=\coprod_{k=0}^{\min(i,j)} S_{j-k,i-k}(T,T')\\
\widetilde S'_{j,i}(T,T')=\coprod_{k=0}^{\min(1,i,j)} S'_{j-k,i-k}(T,T'),
\end{gather*}
where $\coprod$ means the disjoint union.
\end{definition}

\begin{remark}\label{rewriterem}
By definition,
\begin{gather*}
\Braket{D_j U_i T, T'} = \numof{N_{i,j}(T,T')},\\
\Braket{D'_j U_i T, T'} = \numof{N'_{i,j}(T,T')},\\
\Braket{U_j D_i T, T'} = \numof{S'_{i,j}(T,T')},\\
\intertext{and}
\Braket{U_j D'_i T, T'} = \numof{S'_{i,j}(T,T')}
\end{gather*}
for each $T, T'\in\Trees$. Hence 
the equations $(\ref{opeq1})$ and $(\ref{opeq2})$ are respectively equivalent to
the equations
\begin{align*}
\numof{N_{i,j}(T,T')}&=\numof{\widetilde S_{j,i}(T,T')}
\intertext{and}
\numof{N'_{i,j}(T,T')}&=\numof{\widetilde S'_{j,i}(T,T')}.
\end{align*}
\end{remark}

\begin{lemma}\label{rcorreslemma}
For each $T$, $T'\in\Trees$ and each $i$, $j\in\NN$,
there exists bijection from ${N_{i,j}(T,T')}$ to $\widetilde S_{j,i}(T,T')$.
\end{lemma}

\begin{proof}
First we construct 
an element of $\widetilde S_{j,i}(T,T')$
from %$((T,T''),(T',T'')) \in {N_{i,j}(T,T')}$.
an element of ${N_{i,j}(T,T')}$.

Let $((T,T''),(T',T''))$ be an element of $ {N_{i,j}(T,T')}$.
Equivalently, $(T,T'')$ is an edge of $G_{U_i}$ such that 
$T'' \stripminus r_{T'',j} = T'$.
Let $k$ be
$j-\numof{r_{T'',j} \cap  r_{T}}$.
We have $r_{T, j-k}=r_{T'',j} \cap r_{T}$
since $r_{T''}$ is one of the following:
\begin{align*}
 r_{T} ,&\\
 r_{T,l} & \cup \Set{w_{T,l+1}21^{i}|i\leq n},\\
 r_{T,l} & \cup \Set{w_{T,l+1}1^{i}|i\leq n}
\end{align*}
for some $l$, $n\in \NN$.
We consider 
\begin{gather*}
((T\stripminus r_{T, j-k},T),( T\stripminus r_{T, j-k},T')).
\end{gather*}
%We show that 
%$((T,T\stripminus r_{T, k}),(T', T\stripminus r_{T, k})) \in \widetilde S_{j,i}(T,T')$.
%Let $w \in T'' \stripminus r_{T'',j} \setminus T \stripminus r_{T,k}$.
%Since $\nu_{T'',j}(w) \not\in r_{T'',j}$, 
%$\nu_{T'',j}(w) \in T'' \setminus r_{T,k}$.
%Since $w\not\in T\stripminus r_{T,k}$ and $\nu\not\in r_{T,k}$,
%$\nu(w) \in T \setminus T''$.
%Hence $\nu(w)=v1^n$ or $v2^n$ for some $v\in T$.
%Let $v''$ be $\max \Set{u\leq v| u\not\in r_{T,k}}$.
%It follows by definition that
%all nodes $u>v''$ such that $u\leq v$ are in $r_{t,k}\subset r_{T'',k}$.
%Since $v''\not\in r_{T,k} \subset r_{T'',j}$,
%$\nu_{T,k}^{-1}(\Set{v''})\neq \emptyset$ and
%a node $v'\in\nu_{T,k}^{-1}(\Set{v''})$ is also in 
%$\nu_{T'',j}^{-1}(\Set{v''})$.
%Hence, for $v'\in\nu_{T,k}^{-1}(\Set{v''})$, we have 
%$w=v'21^n$ or $v'1^n$ for some $n$.
%Thus, for each
%$w \in T'' \stripminus r_{T'',j} \setminus T \stripminus r_{T,k}$,
%we have $w=v'21^n$ or $v'1^n$ for some $n$ and some $v'T \stripminus r_{T,k}$.
%This means $((T,T\stripminus r_{T, k}),(T', T\stripminus r_{T, k})) \in \widetilde S_{j,i}(T,T')$.
By definition,
$(T\stripminus r_{T, j-k},T)$ is in $G_{D_{j-k}}$.
Since $(T,T'')$ is in $G_{U_i}$,
$( T\stripminus r_{T, j-k},T'' \stripminus r_{T'',j})$ is in $G_{U_{i-k}}$.
Since $T'=T'' \stripminus r_{T'',j}$,
$( T\stripminus r_{T, j-k},T')$ is in $G_{U_{i-k}}$.
Hence we have $((T,T\stripminus r_{T, k}),(T', T\stripminus r_{T, k})) \in \widetilde S_{j,i}(T,T')$.

Next 
 we construct 
an element of ${N_{i,j}(T,T')}$
from %$((T,T''),(T',T'')) \in {N_{i,j}(T,T')}$.
an element of $\widetilde S_{j,i}(T,T')$.
Let $((T''',T),(T''',T'))$ be an element of $\widetilde S_{j,i}(T,T')$.
Equivalently, $(T''',T')$ is an edge of $G_{U_i-k}$ such that 
$T \stripminus r_{T,j-k} = T'''$.

First we consider the case where
$\numof{r_{T}}>j-k$.
Let $\omega = w_{T,j-k+1}$ and 
 $\omega' \in \nu_{T,j-k}^{-1}(\omega)$.
Since $\omega'\in T\stripminus r_{T,j-k}$ and 
$\omega'2 \not\in T\stripminus r_{T,j-k}$,
we have
\begin{gather*}
T'_{\omega'2}=
\Set{\omega'2, \omega'21, \ldots,\omega'21^{n-1}}
\end{gather*}
for some $n\in \NN$, where this is
empty for $n=0$.
For such $n$, let $R$ denote
\begin{gather*}
\Set{\omega 2, \omega 21, \ldots,\omega 21^{n-1+k}},
\end{gather*} where this is
empty for $n+k=0$.

We define $T''$  to be
\begin{gather*}
\widetilde\nu_{T,j-k}(T)\cup r_{T,j-k}\cup R.
\end{gather*}
Since $r_{T''}=r_{T,j-k}\cup R$,
$((T,T''),(T',T''))$ is in 
 ${N_{i,j}(T,T')}$.

Next we consider the case where
$\numof{r_{T}}=j-k$.
Let $\omega$ be 
\begin{gather*}
\max \Set{w\not\in r_{T}|w<w_{T,j-k}}
\end{gather*}
 and
$\omega'\in \nu_{T,j-k}^{-1}(r)$.
Since $\omega'\in T\stripminus r_{T,j-k}$ and 
$\omega'2 \not\in T\stripminus r_{T,j-k}$,
we have
\begin{gather*}
T'_{\omega'2}=
\Set{\omega'2, \omega'21, \ldots,\omega'21^{n-1}}
\end{gather*}
for some $n\in \NN$, where this is
empty  for $n=0$.
For such $n$, let $R$ denote
\begin{gather*}
\Set{w_{T,j-k+1}1, \ldots, w_{T,j-k+1}1^{n-1+k}},
\end{gather*} where this is
empty for $n+k=0$.

We define $T''$  to be
\begin{gather*}
\widetilde\nu_{T,j-k}(T)\cup r_{T,j-k}\cup R.
\end{gather*}
Since $r_{T''}=r_{T,j-k}\cup R$,
$((T,T''),(T',T''))$ is in 
 ${N_{i,j}(T,T')}$.

Thus  we can construct 
an element of ${N_{i,j}(T,T')}$
from %$((T,T''),(T',T'')) \in {N_{i,j}(T,T')}$.
an element of $\widetilde S_{j,i}(T,T')$.

By the definition of them,
 these constructions
are the inverses of each other.
Hence we have the lemma.
\end{proof}

\begin{lemma}\label{dualrcorreslemma}
For $T$, $T'\in\Trees$,
there exists a bijection from $N'_{i,j}(T,T')$ to $\widetilde S'_{j,i}(T,T')$.
\end{lemma}

\begin{proof}
First we construct 
an element of $\widetilde S'_{j,i}(T,T')$
from %$((T,T''),(T',T'')) \in {N_{i,j}(T,T')}$.
an element of ${N'_{i,j}(T,T')}$.

Let $((T,T''),(T',T''))$ be an element of $ {N'_{i,j}(T,T')}$.
Equivalently, $(T,T'')$ is an edge of $G_{U'_i}$ such that 
$T'' \stripminus r_{T'',j} = T'$.
Let $k$ be
$\numof{r_{T'',j} \cap  r_{T}}$.
We have $r_{T, j-k}=r_{T'',j} \cap r_{T}$
since $r_{T''}$ is one of the following:
\begin{align*}
 r_{T} ,&\\
 r_{T,l} & \cup \Set{w_{T,l+1}12^{n}},\\
 r_{T,l} & \cup \Set{w_{T,l+1}2^{n}}
\end{align*}
for some $l$, $n\in \NN$.
It also follows that $k=0$ or $1$.
%It is clear that $r_{T, k}=r_{T'',j} \cap  r_{T}$.
We consider 
\begin{gather*}
((T,T\stripminus r_{T, k}),(T', T\stripminus r_{T, k})).
\end{gather*}
%We show that 
%$((T,T\stripminus r_{T, k}),(T', T\stripminus r_{T, k})) \in \widetilde S'_{j,i}(T,T')$.
%Let $w \in T'' \stripminus r_{T'',j} \setminus T \stripminus r_{T,k}$.
%Since $\nu_{T'',j}(w) \not\in r_{T'',j}$, 
%$\nu_{T'',j}(w) \in T'' \setminus r_{T,k}$.
%Since $w\not\in T\stripminus r_{T,k}$ and $\nu\not\in r_{T,k}$,
%$\nu(w) \in T \setminus T''$.
%Hence $\nu(w)=v12^n$ or $v2^n$ for some $v\in T$.
%Let $v''$ be $\max \Set{u\leq v| u\not\in r_{T,k}}$.
%It follows by definition that
%all nodes $u>v''$ such that $u\leq v$ are in $r_{t,k}\subset r_{T'',k}$.
%Since $v''\not\in r_{T,k} \subset r_{T'',j}$,
%$\nu_{T,k}^{-1}(\Set{v''})\neq \emptyset$ and
%a node $v'\in\nu_{T,k}^{-1}(\Set{v''})$ is also in 
%$\nu_{T'',j}^{-1}(\Set{v''})$.
%Hence, for $v'\in\nu_{T,k}^{-1}(\Set{v''})$, we have 
%$w=v'12^n$ or $v'12^n$ for some $n$.
%Thus, for each
%$w \in T'' \stripminus r_{T'',j} \setminus T \stripminus r_{T,k}$,
%we have $w=v'12^n$ or $v'2^n$ for some $n$ and some $v'T \stripminus r_{T,k}$.
%This means $((T,T\stripminus r_{T, k}),(T', T\stripminus r_{T, k})) \in \widetilde S'_{j,i}(T,T')$.
By definition,
$(T\stripminus r_{T, j-k},T)$ is in $G_{D_{j-k}}$.
Since $(T,T'')$ is in $G_{U'_i}$,
$( T\stripminus r_{T, j-k},T'' \stripminus r_{T'',j})$ is in $G_{U'_{i-k}}$.
Since $T'=T'' \stripminus r_{T'',j}$,
$( T\stripminus r_{T, j-k},T')$ is in $G_{U'_{i-k}}$.
Hence we have $((T,T\stripminus r_{T, k}),(T', T\stripminus r_{T, k})) \in \widetilde S_{j,i}(T,T')$.

Next 
 we construct 
an element of ${N'_{i,j}(T,T')}$
from %$((T,T''),(T',T'')) \in {N_{i,j}(T,T')}$.
an element of $\widetilde S'_{j,i}(T,T')$.
Let $((T''',T),(T''',T'))$ be an element of $\widetilde S'_{j,i}(T,T')$.
Equivalently, $(T''',T')$ is an edge of $G_{U'_i-k}$ such that 
$T \stripminus r_{T,j-k} = T'''$.

First we consider the case where
$\numof{r_{T}}>j-k$.
Let $\omega=w_{T,j-k+1}$
$\omega' \in \nu_{T,j-k}^{-1}(\omega)$.
Since $\omega'\in T\stripminus r_{T,j-k}$ and 
$\omega'2\not\in T\stripminus r_{T,j-k}$,
we have
\begin{gather*}
T'_{\omega'2}=
\Set{\omega'2, \ldots,\omega'2^{n}}
\end{gather*}
for some $n\in \NN$, where this is
empty for $n=0$.
For such $n$, let $R$ denote
\begin{gather*}
\Set{\omega 2, \ldots,\omega 21^{n+k}},
\end{gather*} where this  is
empty for $n+k=0$.

We define $T''$  to be
\begin{gather*}
\widetilde\nu_{T,j-k}(T)\cup r_{T,j-k}\cup R.
\end{gather*}
Since $r_{T''}=r_{T,j-k}\cup \Set{\omega 21^{n+k}}$,
$((T,T''),(T',T''))$ is in 
 ${N'_{i,j}(T,T')}$.

Next we consider the case where
$\numof{r_{T}}=j-k$.
Let $\omega$ be 
\begin{gather*}
\max \Set{w\not\in r_{T}|w<w_{T,j-k}}
\end{gather*}
 and
$\omega'\in \nu_{T,j-k}^{-1}(\omega)$.
Since $\omega'\in T\stripminus r_{T,j-k}$ and 
$\omega'2 \not\in T\stripminus r_{T,j-k}$,
we have
\begin{gather*}
T'_{\omega'2}=
\Set{\omega'2, \ldots,\omega'2^{n}}
\end{gather*}
for some $n\in \NN$, where this is
empty  for $n=0$.
For such $n$, let $R$ denote
\begin{gather*}
\Set{w_{T,j-k}1, w_{T,j-k} 12, \ldots,w_{T,j-k} 12^{n-1+k}},
\end{gather*} where this is
empty  for $n+k=0$.

We define $T''$  to be
\begin{gather*}
\widetilde\nu_{T,j-k}(T)\cup r_{T,j-k}\cup R.
\end{gather*}
Since 
\begin{gather*}
r_{T''}=r_{T,j-k}\cup \Set{w_{T,j-k+1}12^{n-1+k}},
\end{gather*}
$((T,T''),(T',T''))$ is in 
 ${N'_{i,j}(T,T')}$.

Thus  we can construct 
an element of ${N'_{i,j}(T,T')}$
from %$((T,T''),(T',T'')) \in {N_{i,j}(T,T')}$.
an element of $\widetilde S'_{j,i}(T,T')$.

By the definitions of them,
 these constructions
are the inverses of each other.
Hence we have the lemma.
\end{proof}

By Lemmas \ref{rcorreslemma} and \ref{dualrcorreslemma},
we have Theorems \ref{mainthm} and \ref{dualmainthm}.

\section{Application}
In this section, we consider a relation between our generalized Schur
operators and the Loday-Ronco algebra.
\begin{remark}
We have correspondences between
the sets $N_{i,j}(T,T')$ and $\widetilde S_{i,j}(T,T')$
for all $i,j$ by the proof of \ref{rcorreslemma}. 
From them, we can construct
a Robinson-Schensted-Knuth correspondence for paths
of $G_{U}$ and $G_{D}$ by the method in \cite{FomK}.
This correspondence 
is a generalization of the Loday-Ronco correspondence,
which is a Robinson correspondence for binary trees.
%to ``skew'' version.
\end{remark}

\begin{remark}
Maxime Rey  gave a construction of the Loday-Ronco algebra in \cite{Rey}.
He introduced a new  Robinson-Schensted-Knuth correspondence for binary trees 
to construct  the Loday-Ronco algebra.
%We can also construct Robinson-Schensted-Knuth correspondences for paths
%of $G_{U}$ and $G_{D}$ by the method in \cite{FomK}.
%These correspondences are equivalent each other.
Some of our correspondences are equivalent to his correspondence.
\end{remark}

\begin{definition}
For $\lambda,\mu \in V$, we define quasi-symmetric polynomials
$\SD[D]{\lambda}{\mu}(t_1,\ldots,t_n)$,
$\SU[U]{\lambda}{\mu}(t_1,\ldots,t_n)$ and
$\SU[U']{\lambda}{\mu}(t_1,\ldots,t_n)$ by
\begin{align*}
\SD[D]{\lambda}{\mu}(t_1,\ldots,t_n)
&=\Braket{D(t_1)\cdots D(t_n)T,T'}\\
\SU[U]{\lambda}{\mu}(t_1,\ldots,t_n)
&=\Braket{U(t_n)\cdots U(t_1)T',T}\\
\SU[U']{\lambda}{\mu}(t_1,\ldots,t_n)
&=\Braket{U'(t_n)\cdots U'(t_1)T',T}.
\end{align*}
\end{definition}

\begin{remark}\label{rempo}
%Proposition \ref{rempo} follows by the definition of labellings. 
For a labelling $\varphi$ from $T$ to $\{1,\ldots,m\}$,
we define $t^{\varphi}=\prod_{w\in T} t_{\varphi(w)}$.
For a tree $T$, it follows by the definition of labellings that
\begin{align*}
\SD[D]{T}{\emptyset}(t_1,\ldots,t_n)
%\Braket{D(t_n)\cdots D(t_n)T,\emptyset}
&=\sum_{\psi} t^{\psi},\\
\SU[U]{T}{\emptyset}(t_1,\ldots,t_n)
%\Braket{U(t_n)\cdots U(t_1)\emptyset,T}
&=\sum_{\varphi} t^{\varphi},\\
\SU[U']{T}{\emptyset}(t_1,\ldots,t_n)
%\Braket{U'(t_n)\cdots U'(t_1)\emptyset,T}
&=\sum_{\phi} t^{\phi},
\end{align*}
where 
the first sum is over all binary-searching labellings $\psi$ on $T$,
the second sum is over all right-strictly-increasing labellings $\varphi$ on $T$,
and the last sum is over all left-strictly-increasing labellings $\phi$ on $T$.
\end{remark}
\begin{remark}
These polynomials $\SU[U]{T}{\emptyset}(t_1,\ldots,t_n)$ and
$\SD[D]{T}{\emptyset}(t_1,\ldots,t_n)$
are
the commutativizations of the basis elements $\QQQ_{T}$ and $\PPP_{T}$  of $\PPBBTT$ in
Hivert-Novelli-Thibon \cite{hnt}.
\end{remark}

\begin{remark}
Since $D(t)$ and $U(t)$ are generalized Schur operators,
we have 
Pieri's formula
for 
%these polynomials 
$\SU[U]{T}{\emptyset}(t_1,\ldots,t_n)$ and
$\SD[D]{T}{\emptyset}(t_1,\ldots,t_n)$
by \cite{pieri}. 
By \cite{FomK}, 
we have
Cauchy identity
for 
them.
We also have
a ``skew'' version of them.

We also have
Pieri's formula and Cauchy identity
for 
$\SU[U']{T}{\emptyset}(t_1,\ldots,t_n)$ and
$\SD[D]{T}{\emptyset}(t_1,\ldots,t_n)$
\end{remark}

\begin{remark}
These polynomials are not symmetric in general.
For example, since
\begin{align*}
&D(t_1)D(t_2) \{ 0, 1 , 12\}\\
&=D(t_1)( \{ 0, 1 , 12\} + t_2 \{0,2\} + t_2^2 \{0\})\\
&=( \{ 0, 1 , 12\} + t_1 \{0,2\} + t_1^2 \{0\})
+t_2(\{0,2\}+ t_1 \{0\}) 
+t_2^2(\{0\} +t_1\emptyset),
\end{align*}
$\Braket{D(t_n)\cdots D(t_n)T,\emptyset}=t_1 t_2^2$ is not symmetric for
 $T=\{ 0, 1 , 12\}$.
The fact that
$D_i$ does not commute with $D_j$  in general
implies this fact.
\end{remark}


\begin{thebibliography}{99}
\bibitem{FomRSK}
S. Fomin,
\emph{Generalized Robinson-Schensted-Knuth correspondence},
Zariski Nauchn. Sem. Leningrad. Otdel. Mat. Inst. Steklov. (LOMI)
\textbf{155}
(1986),
156--175, 195 (Russian);
English transl., 
J. Soviet Math.
41(1988),
979--991.

\bibitem{FomD}
S. Fomin, 
\emph{Duality of graded graphs}, 
J. Algebraic Combin. 
\textbf{3}
(1994),
357--404.

\bibitem{FomS}
S. Fomin,
\emph{Schensted algorithms for dual graded graphs},
J. Algebraic Combin.
\textbf{4}
(1995),
5--45.


\bibitem{FomK}
S. Fomin, 
\emph{Schur operators and Knuth correspondences}, 
J. Combin. Theory,
Ser. A 
\textbf{72} 
(1995), 
277--292.
%\bibitem{Ful}
%Fulton,W.,
%Young Tableaux; with applications to representation theory and geometry,
%volume 35 of London Mathematical Society Student Texts,
%Cambridge University Press,
%New York,
%1997.

\bibitem{Gess}
Ira M. Gessel, 
\emph{Counting paths in Young's lattice},
J. Statistical planning and inference.
\textbf{34}
(1993),
125--134.


\bibitem{hnt}
F. Hivert, J. Novelli and J. Thibon, 
\emph{The algebra of binary search trees},
Theor. Comput. Sci. 339, 1 (Jun. 2005), 129--165. 
DOI=\url{http://dx.doi.org/10.1016/j.tcs.2005.01.012}

\bibitem{pieri}
Y. Numata, 
\emph{Pieri's Formula for Generalized Schur Polynomials},
preprint,
\url{arXiv:math.CO/0606386}.

\bibitem{Nzc}
J. Nzeutchap,
Graded Graphs and Fomin's $r$-correspondences associated to the Hopf
	Algebras of Planar Binary Trees, Quasi-symmetric Functions and
	Noncommutative Symmetric Functions,
FPSAC '06,
2006.
\url{http://garsia.math.yorku.ca/fpsac06/papers/53.pdf}


\bibitem{Rey}
M. Rey,
A new construction of the Loday-Ronco algebra,
FPSAC '06,
2006.
\url{http://garsia.math.yorku.ca/fpsac06/papers/51.pdf}


\bibitem{Rob}
T. Roby,
\emph{Applications and extensions of Fomin's generalization of the Robinson-Schensted 
 correspondence to differential posets},
Ph.D. thesis,
M.I.T.,
1991.


\bibitem{StaD}
R. Stanley,
\emph{Differential posets},
J. American Math. Soc,
\textbf{1}
(1988),
919--961.

\bibitem{StaV}
R. Stanley,
\emph{Variations on differential posets},
Invariant theory and tableaux
(Stanton,D.,ed.),
IMA volumes in mathematics
and its applications,
Springer-Verlag,
New York,
145--165.

\end{thebibliography}
\end{document}